# ALGORITHMS FOR SOLVING LINEAR CONGRUENCES AND SYSTEMS OF LINEAR CONGRUENCES

Florentin Smarandache
University of New Mexico
200 College Road
Gallup, NM 87301, USA
E-mail: smarand@unm.edu

In this article we determine several theorems and methods for solving linear congruences and systems of linear congruences and we find the number of distinct solutions. Many examples of solving congruences are given.

**§1. Properties for solving linear congruences.**

**Theorem 1.** The linear congruence $a_1 x_1 + ... + a_n x_n \equiv b \pmod{m}$ has solutions if and only if $(a_1,...,a_n,m) | b$.

*Proof:*
$a_1 x_1 + ... + a_n x_n \equiv b \pmod{m} \Leftrightarrow a_1 x_1 + ... + a_n x_n - my = b$ is a linear equation which has solutions in the set of integer numbers $\Leftrightarrow (a_1,...,a_n,-m) | b \Leftrightarrow (a_1,...,a_n,m) | b$.

If $m = 0$, $a_1 x_1 + ... + a_n x_n \equiv b \pmod{0} \Leftrightarrow a_1 x_1 + ... + a_n x_n = b$ has solutions in the set of integer numbers $\Leftrightarrow (a_1,...,a_n) | b \Leftrightarrow (a_1,...,a_n,0) | b$.

**Theorem 2.** The congruence $ax \equiv b \pmod{m}$, $m \neq 0$, with $(a,m) = d | b$, has $d$ distinct solutions.

The proof is different of that from the number's theory courses:
$ax \equiv b \pmod{m} \Leftrightarrow ax - my = b$ has solutions in the set of integer numbers; because $(a,m) = d | b$ it results: $a = a_1 d$, $m = m_1 d$, $b = b_1 d$ and $(a_1, m_1) = 1$, $a_1 dx - m_1 dy = b_1 d \Leftrightarrow a_1 x - m_1 y = b_1$. Because $(a_1, m_1) = 1$ it results that the general solution of this equation is $\begin{cases} x = m_1 k_1 + x_0 \\ y = a_1 k_1 + y_0 \end{cases}$, where $k_1$ is a parameter and $k_1 \in \mathbb{Z}$, and where $(x_0, y_0)$ constitutes a particular solution in the set of integer numbers of this equation; $x = m_1 k_1 + x_0$, $k_1 \in \mathbb{Z}$, $m_1, x_0 \in \mathbb{Z} \Rightarrow x \equiv m_1 k_1 + x_0 \pmod{m}$. We'll assign values to $k_1$ to find all the solutions of the congruence.
It is evident that $k_1 \in \{0,1,2,...,d-1,d,d+1,...,m-1\}$ which constitutes a complete system of residues modulo $m$.
(Because $ax \equiv b \pmod{m} \Leftrightarrow ax \equiv b \pmod{-m}$, we suppose $m > 0$.)

Let $D = \{0,1,2,...,d-1\}$; $D \subseteq M$, $\forall \alpha \in M$, $\exists \beta \in D : \alpha \equiv \beta \pmod{d} | m_1$
(because $D$ constitutes a complete system of residues modulo $d$).

It results that $\alpha m_1 \equiv \beta m_1 \pmod{dm_1}$; because $x_0 \equiv x_0 \pmod{dm_1}$, it results:
$$m_1 \alpha + x_0 \equiv m_1 \beta + x_0 \pmod{m}.$$

Therefore $\forall \alpha \in M$, $\exists \beta \in D: m_1\alpha + x_0 \equiv m_1\beta + x_0 (\mod m)$; thus $k_1 \in D$.
$\forall \gamma, \delta \in D$, $\gamma \not\equiv \delta (\mod d) \mid m_1 \Rightarrow \gamma m_1 \not\equiv \delta m_1 (\mod dm_1)$; $m_1 \neq 0$. It results that $m_1\gamma + x_0 \equiv m_1\delta + x_0 (\mod m)$ is false, that is, we have exactly $cardD = d$ distinct solutions.

**Remark 1**. If $m = 0$, the congruence $ax \equiv b(\mod 0)$ has one solution if $a \mid b$; otherwise it does not have solutions.

*Proof:*
$ax \equiv b(\mod 0) \Leftrightarrow ax = b$ has a solution in the set of integer numbers $\Leftrightarrow a \mid b$.

**Theorem 3.** (A generalization of the previous theorem)
The congruence $a_1x_1 + ... + a_nx_n \equiv b(\mod m)$, $m_1 \neq 0$, with $(a_1,...,a_n,m) = d \mid b$ has $d \cdot |m|^{n-1}$ distinct solutions.

*Proof:*
Because $a_1x_1 + ... + a_nx_n \equiv b(\mod m) \Leftrightarrow a_1x_1 + ... + a_nx_n \equiv b(\mod -m)$, we can consider $m > 0$.

The proof is done by induction on $n =$ the number of variables.
For $n = 1$ the affirmation is true in conformity with theorem 2.
Suppose that it is true for $n - 1$. Let's proof that it is true for $n$.
Let the congruence with $n$ variables $a_1x_1 + ... + a_nx_n \equiv b(\mod m)$, $a_1x_1 + ... + a_{n-1}x_{n-1} \equiv b - a_nx_n (\mod m)$. If we consider that $x_n$ is fixed, the congruence $a_1x_1 + ... + a_{n-1}x_{n-1} \equiv b - a_nx_n (\mod m)$ is a congruence with $n - 1$ variables. To have solutions we must have $(a_1,...,a_{n-1},m) = \delta \mid b - a_nx_n \Leftrightarrow b - a_nx_n \equiv 0(\mod \delta)$.

Because $\delta \mid m \Rightarrow \frac{m}{\delta} \in \mathbb{Z}$, therefore we can multiply the previous congruence with $\frac{m}{\delta}$. It results that

$$\frac{ma_n}{\delta}x_n \equiv \frac{mb}{\delta}(\mod \delta \cdot \frac{m}{\delta}) \quad (*)$$

which has $\left(\frac{ma_n}{\delta}, \delta\frac{m}{\delta}\right) = \frac{m}{\delta}(a_n, \delta) = \frac{m}{\delta}(a_n,(a_1,...,a_{n-1},m)) = \frac{m}{\delta}(a_1,...,a_{n-1},a_n,m)\frac{m}{\delta} \cdot d$ distinct solutions for $x_n$. Let $x_n^0$ be a particular solution of the congruence (*). It results that $a_1x_1 + ... + a_{n-1}x_{n-1} \equiv b - a_nx_n^0 (\mod m)$ has, conform to the induction's hypothesis, $\delta \cdot m^{n-2}$ distinct solutions for $x_1,...,x_{n-1}$ where $\delta = (a_1,...,a_{n-1},m)$.

Therefore the congruence $a_1x_1 + ... + a_{n-1}x_{n-1} + a_nx_n \equiv b(\mod m)$ has $\frac{m}{\delta} \cdot d \cdot \delta \cdot m^{n-2} = d \cdot m^{n-1}$ distinct solutions for $x_1,...,x_{n-1}$ and $x_n$.

## §2. A METHOD FOR SOLVING LINEAR CONGRUENCES

Let's consider the congruence $a_1x_1 + ... + a_nx_n \equiv b(\mod m)$, $m \neq 0$,

$a_i \equiv a_i'(\mod m)$ and $b \equiv b'(\mod m)$ with $0 \le a_i'$, $b \le m-1$ (we made the nonrestrictive hypothesis $m > 0$). We obtain:

$a_1 x_1 + ... + a_n x_n \equiv b(\mod m) \Leftrightarrow a_1' x_1 + ... + a_n' x_n \equiv b'(\mod m)$, which is a linear equation; when it is resolved in $\mathbb{Z}$ it has the general solution:

$$\begin{cases} x_1 = \alpha_{11} k_1 + ... + \alpha_{1n} k_n + \gamma_1 \\ \vdots \\ x_n = \alpha_{n1} k_1 + ... + \alpha_{nn} k_n + \gamma_n \\ y = \alpha_{n+1,1} k_1 + ... + \alpha_{n+1,n} k_n + \gamma_{n+1} \end{cases}$$

$k_j$ being parameters $\in \mathbb{Z}$, $j = \overline{1,n}$, $\alpha_{ij}, \gamma_i \in \mathbb{Z}$, constants, $i = \overline{1,n+1}$, $j = \overline{1,n}$.

Let's consider $\alpha_{ij}' \equiv \alpha_{ij}(\mod m)$ and $\gamma_i' \equiv \gamma_i(\mod m)$ with $0 \le \alpha_{ij}'$, $\gamma' \le m-1$; $i = \overline{1, n+1}$, $j = \overline{1,n}$.

Therefore

$$\begin{cases} x_1 = \alpha_{11}' k_1 + ... + \alpha_{1n}' k_n + \gamma_1' (\mod m) \\ \vdots \\ x_n = \alpha_{n1}' k_1 + ... + \alpha_{nn}' k_n + \gamma_n' (\mod m) \end{cases} \quad ; k_j = \text{parameters} \in \mathbb{Z}, j = \overline{1,n}; \quad (**)$$

Let's consider $(\alpha_{1j}', ..., \alpha_{nj}', m) = d_j$, $j \in \overline{1,n}$. We'll prove that for $k_j$ it would be sufficient to only give the values $0, 1, 2, ..., \dfrac{m}{d_j} - 1$; for $k_j = \dfrac{m}{d_j} - 1 + \beta'$ with $\beta' \ge 1$ we obtain $k_j = \dfrac{m}{d_j} + \beta$ with $\beta \ge 0$; $\beta', \beta \in \mathbb{Z}$.

$\alpha_{ij}' k_j = \alpha_{ij}'' d_j k_j = \alpha_{ij}'' m + \alpha_{ij}'' d_j \beta \equiv \alpha_{ij}'' d_j \beta (\mod m)$; we denoted $\alpha_{ij}' = \alpha_{ij}'' d_j$ because $d_j | \alpha_{ij}'$.

We make the notation $m = d_j m_j$, $m_j = \dfrac{m}{d_j}$.

Let's consider $\eta \in \mathbb{Z}$, $0 \le \eta \le m-1$ such that $\eta = \alpha_{ij}'' d_j \beta(\mod d_j m_j)$; it results $d_j | \eta$.

Therefore $\eta = d_j \gamma$ with $0 \le \gamma \le m_{j-1}$ because we have that $d_j \gamma \equiv \alpha_{ij}'' d_j(\mod d_j m_j)$, which is equivalent to $\gamma \equiv \alpha_{ij}'' \beta(\mod m_j)$.

Therefore $\forall k_j \in \mathbb{N}$, $\exists \gamma \in \{0, 1, 2, ..., m_{j-1}\}$: $\alpha_{ij}' k_j \equiv d_j \gamma(\mod m)$; analogously, if the parameter $k_j \in \mathbb{Z}$. Therefore $k_j$ takes values from $0, 1, 2, ...$ to at most $m_j - 1$; $j \in \overline{1,n}$.

Through this parameterization for each $k_j$ in (**), we obtain the solutions of the linear congruence. We eliminate the repetitive solutions. We obtain exactly $d \cdot |m|^{n-1}$ distinct solutions.

**Example 1.** Let's resolve the following linear congruence:

$$2x + 7y - 6z \equiv -3 \pmod{4}$$

*Solution:* $7 \equiv 3 \pmod 4$, $-6 \equiv 2 \pmod 4$, $-3 \equiv 1 \pmod 4$.

It results that $2x + 3y + 2z \equiv 1 \pmod 4$; $(2,3,2,4) = 1 | 1$ therefore the congruence has solutions and it has $1 \cdot 4^{3-1} = 16$ distinct solutions.

The equation $2x + 3y + 2z - 4t = 1$ resolved in integer numbers, has the general solution:
$$\begin{cases} x = 3k_1 - k_2 - 2k_3 - 1 \equiv 3k_1 + 3k_2 + 2k_3 + 3 \pmod 4 \\ y = -2k_1 \phantom{aaaaaaa} + 1 \equiv 2k_1 \phantom{aaaaaaaa} + 1 \pmod 4 \\ z = \phantom{aaaa} k_2 \phantom{aaaaaaa} \equiv \phantom{aaa} k_2 \phantom{aaaaaa} \pmod 4 \end{cases}$$

$k_j$ are parameters $\in \mathbb{Z}$, $j = \overline{1,3}$.

(We did not write the expression for $t$, because it doesn't interest us).

We assign values to the parameters. $k_j$ takes values from 0 to at most $m_j - 1$;

$k_3$ takes values from 0 to $m_3 - 1 = \dfrac{m}{d_3} - 1 = \dfrac{4}{(2,0,0)} - 1 = \dfrac{4}{2} - 1 = 1$;

$$k_3 = 0 \Rightarrow \begin{pmatrix} x \equiv 3k_1 + 3k_2 + 3 \pmod 4 \\ y \equiv 2k_1 \phantom{aaaaa} + 1 \pmod 4 \\ z \equiv \phantom{aaa} k_2 \phantom{aaaa} \pmod 4 \end{pmatrix};$$

$$k_3 = 1 \Rightarrow \begin{pmatrix} 3k_1 + 3k_2 + 1 \\ 2k_1 \phantom{aaaa} + 1 \\ \phantom{aa} k_2 \end{pmatrix}$$

$k_1$ takes values from 0 to at most 3.

$$k_1 = 0 \Rightarrow \begin{pmatrix} 3k_2 + 3 \\ 1 \\ k_2 \end{pmatrix}, \begin{pmatrix} 3k_2 + 1 \\ 1 \\ k_2 \end{pmatrix}; \quad k_1 = 1 \Rightarrow \begin{pmatrix} 3k_2 + 2 \\ 3 \\ k_2 \end{pmatrix}, \begin{pmatrix} 3k_2 \\ 3 \\ k_2 \end{pmatrix};$$

for $k_1 = 2$ and 3 we obtain the same expressions as for $k_1 = 1$ and 0.

$k_2$ takes values from 0 to at most 3.

$$k_2 = 0 \Rightarrow \begin{pmatrix} 3 \\ 1 \\ 0 \end{pmatrix}, \begin{pmatrix} 1 \\ 1 \\ 0 \end{pmatrix}, \begin{pmatrix} 2 \\ 3 \\ 0 \end{pmatrix}, \begin{pmatrix} 0 \\ 3 \\ 0 \end{pmatrix}; \quad k_2 = 2 \Rightarrow \begin{pmatrix} 1 \\ 1 \\ 2 \end{pmatrix}, \begin{pmatrix} 3 \\ 1 \\ 2 \end{pmatrix}, \begin{pmatrix} 0 \\ 3 \\ 2 \end{pmatrix}, \begin{pmatrix} 2 \\ 3 \\ 2 \end{pmatrix};$$

$$k_2 = 1 \Rightarrow \begin{pmatrix} 2 \\ 1 \\ 1 \end{pmatrix}, \begin{pmatrix} 0 \\ 1 \\ 1 \end{pmatrix}, \begin{pmatrix} 1 \\ 3 \\ 1 \end{pmatrix}, \begin{pmatrix} 3 \\ 3 \\ 1 \end{pmatrix}; \quad k_2 = 3 \Rightarrow \begin{pmatrix} 0 \\ 1 \\ 3 \end{pmatrix}, \begin{pmatrix} 2 \\ 1 \\ 3 \end{pmatrix}, \begin{pmatrix} 3 \\ 3 \\ 3 \end{pmatrix}, \begin{pmatrix} 1 \\ 3 \\ 3 \end{pmatrix};$$

which represent all distinct solutions of the congruence.

**Remark 2.** By simplification or amplification of the congruence (the division or multiplication with a number $\neq 0, 1, -1$), which affects also the module, we lose solutions, respectively foreign solutions are introduced.

**Example 2.**

1) The congruence $2x - 2y \equiv 6 \pmod{4}$ has the solutions

$$\binom{3}{0}, \binom{1}{0}, \binom{0}{1}, \binom{2}{1}, \binom{1}{2}, \binom{3}{2}, \binom{2}{3}, \binom{0}{3};$$

2) If we would simplify by 2, we would obtain the congruence $x - y \equiv 3 \pmod{2}$, which has the solutions $\binom{1}{0}, \binom{0}{1}$; therefore we lose solutions.

3) If we would amplify with 2, we would obtain the congruence $4x - 4y \equiv 12 \pmod{4}$, which has the solutions:

$$\binom{3}{0}, \binom{5}{0}, \binom{7}{0}, \binom{1}{0}, \binom{4}{1}, \binom{6}{1}, \binom{0}{1}, \binom{2}{1},$$

$$\binom{5}{2}, \binom{7}{2}, \binom{1}{2}, \binom{3}{2}, \binom{6}{3}, \binom{0}{3}, \binom{2}{3}, \binom{4}{3},$$

$$\binom{7}{4}, \binom{1}{4}, \binom{3}{4}, \binom{5}{4}, \binom{0}{5}, \binom{2}{5}, \binom{4}{5}, \binom{6}{5},$$

$$\binom{1}{6}, \binom{3}{6}, \binom{5}{6}, \binom{7}{6}, \binom{2}{7}, \binom{4}{7}, \binom{6}{7}, \binom{0}{7};$$

therefore we introduce foreign solutions.

**Remark 3.** By the division or multiplication of a congruence with a number which is prime with the module, without dividing or multiplying the module, we obtain a congruence which has the same solutions with the initial one.

**Example 3.** The congruence $2x + 3y \equiv 2 \pmod{5}$ has the same solutions as the congruence $6x + 9y \equiv 6 \pmod{5}$ as follows:

$$\binom{0}{1}, \binom{2}{1}, \binom{3}{2}, \binom{4}{3}, \binom{0}{4}.$$

## §2. PROPERTIES FOR SOLVING SYSTEMS OF LINEAR CONGRUENCES.

In this paragraph we will obtain some interesting theorems regarding the systems of congruences and then a method of solving them.

**Theorem 1.** The system of linear congruences:

(1) $a_{i1}x_1 + ... + a_{in}x_n \equiv b \pmod{m_i}$, $i = \overline{1,r}$, has solutions if and only if the system of linear equations:

(2) $a_{i1}x_1 + \ldots + a_{in}x_n - m_i y_i = b$, $y_i$ unknowns $\in \mathbb{Z}$, $i = \overline{1,r}$, has solutions in the set of integer numbers.

The proof is evident.

**Remark 1.** From the anterior theorem it results that to solve the system of congruences (1) is equivalent with solving in integer numbers the system of linear equations (2).

**Theorem 2.** (A generalization of the theorem from p. 20, from [1]).

The system of congruences $a_i x \equiv b_i (\bmod\, m_i)$, $m_i \neq 0$, $i = \overline{1,r}$ admits solutions if and only if: $(a_i, m_i) | b_i$, $i = \overline{1,r}$ and $(a_i m_j, a_j m_i)$ divides $a_i b_j - a_j b_i$, $i, j = \overline{1,r}$.

*Poof:*

$\forall i = \overline{1,r}$, $a_i x \equiv b_i (\bmod\, m_i)$ $\Leftrightarrow$ $\forall i = \overline{1,r}$, $a_i x = b_i + m_i y_i$, $y_i$ being unknowns $\in \mathbb{Z}$; these Diophantine equations, taken separately, have solutions if and only if $(a_i, m_i) | b_i$, $i = \overline{1,r}$.

$\forall i, j = \overline{1,r}$, from: $a_i x = b_i + y_i m_i\, |\, a_j$ and $a_j \cdot x = b_j + y_j \cdot m_j\, |\, a_i$ we obtain: $a_i a_j \cdot x = a_j b_i + a_j \cdot m_i y_i = a_i b_j + a_i \cdot m_j y_j$, Diophantine equations which have solution if and only if $(a_i m_j, a_j m_i) | a_i b_j - a_j b_i$, $i, j = \overline{1,r}$.

**Consequence.** (We obtain a simpler form for the theorem from p. 20 of [1]). The system of congruences $x \equiv b_i (\bmod\, m_i)$, $m_i \neq 0$, $i = \overline{1,r}$ has solutions if and only if $(m_i, m_j) | b_i - b_j$, $i, j = \overline{1,r}$.

*Proof:*

From theorem 2, $a_i = 1$, $\forall i = \overline{1,r}$ and $(1, m_i) = 1 | b_i$, $i = \overline{1,r}$.

## §4. METHOD FOR SOLVING SYSTEMS OF LINEAR CONGRUENCES

Let's consider the system of linear congruences:

(3) $a_{i1}x_1 + a_{i2}x_2 + \ldots + a_{in} \equiv b_i\, (\bmod\, m_i)$, $i = \overline{1,r}$, the system's matrix rank being $r < n$, $a_{ij}$, $b_i$, $m_i \in \mathbb{Z}$, $m_i \neq 0$, $i = \overline{1,r}$, $j = \overline{1,n}$.

According to §1 from this chapter, we can consider:

(*) $0 \leq a_{ij} \leq |m_i| - 1$, $0 \leq b_i \leq |m_i| - 1$, $\forall i = \overline{1,r}$, $j = \overline{1,n}$. From the theorem 1 and the remark 1 it results that, to solve this system of congruences is equivalent with solving in integer numbers the system of equations:

(4) $a_{i1}x_1 + \ldots + a_{in}x_n - m_i y_i = b_i$, $i = \overline{1,r}$, the system's matrix rank being $r < n$.

Using the algorithm from [2], we obtain the general solution of this system:

$$\begin{cases} x_1 = \alpha_{11}k_1 + \ldots + \alpha_{1n}k_n + \beta_1 \\ \ldots\ldots\ldots\ldots\ldots\ldots\ldots\ldots\ldots\ldots \\ x_n = \alpha_{n1}k_1 + \ldots + \alpha_{nn}k_n + \beta_n \\ y_1 = \alpha_{n+1,1}k_1 + \ldots + \alpha_{n+1,n}k_n + \beta_{n+1} \\ \ldots\ldots\ldots\ldots\ldots\ldots\ldots\ldots\ldots\ldots \\ y_r = \alpha_{n+r,1}k_1 + \ldots + \alpha_{n+r,n}k_n + \beta_{n+r} \end{cases}$$

$\alpha_{hj}, \beta_h \in \mathbb{Z}$ and $k_j$ are parameters $\in \mathbb{Z}$.

Let's consider $m = [m_1, \ldots, m_r] > 0$; because the variables $y_1, \ldots, y_r$ don't interest us, we'll retain only the expressions of $x_1, \ldots, x_n$.
Therefore:

(5) $x_i = \alpha_{i1}k_1 + \ldots + \alpha_{in}k_n + \beta_i$, $i = \overline{1,n}$ and again we can suppose that

(**) $0 \le \alpha_{hj} \le m-1$, $0 \le \beta_h \le m-1$, $h = \overline{1,n}$, $j = \overline{1,n}$.

We have: $x_i \equiv \alpha_{i1}k_1 + \ldots + \alpha_{in}k_n + \beta_i \pmod{m}$, $i = \overline{1,n}$. Evidently $k_j$ takes the values of at most the integer numbers from 0 to $m-1$. Conform to the same observations from §1 from this chapter, for $k_j$ it is sufficient to give only the values $0, 1, 2, \ldots, \dfrac{m}{d_j} - 1$

where

(***) $d_j = \left(\alpha_{1j}, \ldots, \alpha_{nj}, m\right)$, for any $j = \overline{1,n}$.

By the parameterization of $k_1, \ldots, k_n$ in (5) we obtain all the solutions of the system of linear congruence (1); $k_j$ takes at most the values $0, 1, 2, \ldots, \dfrac{m}{d_j} - 1$; we eliminate the repeating solutions.

**Remark 2.** The considerations (\*), (\*\*), and (\*\*\*) have the roll of making the calculation easier, to reduce the computational volume. This algorithm of solving the linear congruence works also without these considerations, but it is more difficult.

**Example.** Let's solve the following system of linear congruences:

(6) $\begin{cases} 3x + 7y - z \equiv 2 \pmod{2} \\ 5y - 2z \equiv 1 \pmod{3} \end{cases}$

*Solution:* The system of linear congruences (6) is equivalent with:

(7) $\begin{cases} x + y + z \equiv 0 \pmod{2} \\ 2y + z \equiv 1 \pmod{3} \end{cases}$

which is equivalent with the system of linear equations:

(8) $\begin{cases} x + y + z - 2t_1 = 0 \\ 2y + z - 3t_2 = 1 \end{cases}$

$x, y, z, t_1, t_2$ unknowns $\in \mathbb{Z}$

This has the general solution (see [2]):

$$\begin{cases} x = -2k_1 + 2k_2 + 3k_3 + 1 \\ y = \phantom{-2}k_1 \phantom{+2k_2} - 3k_3 - 1 \\ z = \phantom{-2}k_1 \\ t_1 = \phantom{-2k_1+}k_2 \\ t_2 = \phantom{-2k_1+2k_2+}k_3 \end{cases}$$

where $k_1$, $k_2$, $k_3$ are parameters $\in \mathbb{Z}$.

The values of $t_1$ and $t_2$ don't interest us; $m = [2, 3] = 6$. Therefore:

$$\begin{cases} x \equiv 4k_1 + 2k_2 + 3k_3 + 1 \pmod 6 \\ y \equiv \phantom{4}k_1 \phantom{+ 2k_2} + 3k_3 + 5 \pmod 6 \\ z \equiv \phantom{4}k_1 \phantom{+ 2k_2 + 3k_3 +} \pmod 6 \end{cases}$$

$k_3$ takes values from 0 to $\dfrac{6}{(3,3,0,6)} - 1 = 1$; $k_2$ from 0 to 2; $k_1$ from 0 to at most 5.

$$k_3 = 0 \Rightarrow \begin{pmatrix} x \equiv 4k_1 + 2k_2 + 1 \pmod 6 \\ y \equiv \phantom{4}k_1 \phantom{+ 2k_2} + 5 \pmod 6 \\ z \equiv \phantom{4}k_1 \phantom{+ 2k_2 + 1} \pmod 6 \end{pmatrix};$$

$$k_3 = 1 \Rightarrow \begin{pmatrix} 4k_1 + 2k_2 + 4 \\ k_1 \phantom{+ 2k_2} + 2 \\ k_1 \end{pmatrix};$$

$$k_2 = 0, 1, 2 \Rightarrow \begin{pmatrix} 4k_1 + 1 \\ k_1 + 5 \\ k_1 \end{pmatrix}, \begin{pmatrix} 4k_1 + 4 \\ k_1 + 2 \\ k_1 \end{pmatrix}, \begin{pmatrix} 4k_1 + 3 \\ k_1 + 5 \\ k_1 \end{pmatrix}, \begin{pmatrix} 4k_1 \\ k_1 + 2 \\ k_1 \end{pmatrix}, \begin{pmatrix} 4k_1 + 5 \\ k_1 + 5 \\ k_1 \end{pmatrix}, \begin{pmatrix} 4k_1 + 2 \\ k_1 + 2 \\ k_1 \end{pmatrix};$$

$k_1 = 0, 1, 2, 3, 4, 5 \Rightarrow$

$$\begin{pmatrix} 1 \\ 5 \\ 0 \end{pmatrix}, \begin{pmatrix} 4 \\ 2 \\ 0 \end{pmatrix}, \begin{pmatrix} 3 \\ 5 \\ 0 \end{pmatrix}, \begin{pmatrix} 0 \\ 2 \\ 0 \end{pmatrix}, \begin{pmatrix} 5 \\ 5 \\ 0 \end{pmatrix}, \begin{pmatrix} 2 \\ 2 \\ 0 \end{pmatrix}, \begin{pmatrix} 5 \\ 0 \\ 1 \end{pmatrix}, \begin{pmatrix} 2 \\ 3 \\ 1 \end{pmatrix}, \begin{pmatrix} 1 \\ 0 \\ 1 \end{pmatrix}, \begin{pmatrix} 4 \\ 3 \\ 1 \end{pmatrix}, \begin{pmatrix} 3 \\ 0 \\ 1 \end{pmatrix}, \begin{pmatrix} 0 \\ 3 \\ 1 \end{pmatrix},$$

$$\begin{pmatrix} 3 \\ 1 \\ 2 \end{pmatrix}, \begin{pmatrix} 0 \\ 4 \\ 2 \end{pmatrix}, \begin{pmatrix} 5 \\ 1 \\ 2 \end{pmatrix}, \begin{pmatrix} 2 \\ 4 \\ 2 \end{pmatrix}, \begin{pmatrix} 1 \\ 1 \\ 2 \end{pmatrix}, \begin{pmatrix} 4 \\ 4 \\ 2 \end{pmatrix}, \begin{pmatrix} 1 \\ 2 \\ 3 \end{pmatrix}, \begin{pmatrix} 4 \\ 5 \\ 3 \end{pmatrix}, \begin{pmatrix} 3 \\ 2 \\ 3 \end{pmatrix}, \begin{pmatrix} 0 \\ 5 \\ 3 \end{pmatrix}, \begin{pmatrix} 5 \\ 2 \\ 3 \end{pmatrix}, \begin{pmatrix} 2 \\ 5 \\ 3 \end{pmatrix},$$

$$\begin{pmatrix} 5 \\ 3 \\ 4 \end{pmatrix}, \begin{pmatrix} 2 \\ 0 \\ 4 \end{pmatrix}, \begin{pmatrix} 1 \\ 3 \\ 4 \end{pmatrix}, \begin{pmatrix} 4 \\ 0 \\ 4 \end{pmatrix}, \begin{pmatrix} 3 \\ 3 \\ 4 \end{pmatrix}, \begin{pmatrix} 0 \\ 0 \\ 4 \end{pmatrix}, \begin{pmatrix} 3 \\ 4 \\ 5 \end{pmatrix}, \begin{pmatrix} 0 \\ 1 \\ 5 \end{pmatrix}, \begin{pmatrix} 5 \\ 4 \\ 5 \end{pmatrix}, \begin{pmatrix} 2 \\ 1 \\ 5 \end{pmatrix}, \begin{pmatrix} 1 \\ 4 \\ 5 \end{pmatrix}, \begin{pmatrix} 4 \\ 1 \\ 5 \end{pmatrix};$$

which constitute the 36 distinct solutions of the system of linear congruences (6).

**REFERENTES:**

[1] Constantin P. Popovici – "Curs de teoria numerelor", EDP, Bucureşti, 1973.
[2] Florentin Smarandache – "Integer algorithms to solve linear equations and systems", Ed. Scientifique, Casablanca, 1984.